\documentclass[10pt,twoside]{article}
\usepackage{hyperref}
\usepackage{amssymb}
\usepackage[mathscr]{eucal}
\usepackage{eufrak}
\usepackage{amsmath}
\usepackage{mathrsfs}

\def\intl#1{\int\limits_{#1}}
\def\intll#1#2{\int\limits_{#1}^{#2}}
\def\dm{|\hskip-0.05cm|}
\def\ssfrac#1#2{\mbox{\large$\frac{#1}{#2}$}}
\def\sfrac#1#2{\mbox{\Large$\frac{#1}{#2}$}}
\def\n{\nabla}
\def\tm{|\hskip-0.05cm|\hskip-0.05cm|}
\def\displ{\displaystyle}
\def\VSE{\vspace{6pt}\\&\displ }
\def\VS{\vspace{6pt}\\\displ }
\def\rf#1{{\rm(\ref{#1})}}
\def\chiu{\hfill$\displaystyle\vspace{4pt}
\underset\Box\null$\par}
\def\P{{\bf Proof. }}
\def\O{\Omega}

\def\R{\Bbb R}
\def\N{\Bbb N}
\def\o{\"{o}}
\def\à{\`{a}}
\def\è{\`{e}}
\def\ì{\`{i}}
\def\ù{\`{u}}
\def\ò{\`{o}}
\def\é{\'{e}}

\def\dy{\displaystyle}
\def\vep{\varepsilon}

\def\be{\begin{equation}}
\def\ba{\begin{array}}
\def\ea{\end{array}}
\def\ee{\end{equation}}
\def\vs1{\vspace{1ex}}

\def\ov{\overline}
\def\po{\partial\Omega}
\font\sc=cmcsc10
\title{A new proof of existence in the $L^3$-setting of solutions to the Navier-Stokes Cauchy problem}
\author{\sc F. Crispo and P. Maremonti
\thanks{
Dipartimento di Matematica e Fisica, 
Universit\`{a} degli Studi della Campania
 ``L. Vanvitelli'', via Vivaldi 43, 81100 Caserta,
 Italy.
francesca.crispo@unicampania.it;
paolo.maremonti@unicampania.it\newline This research was performed under the auspices of the group GNFM-INdAM. The research activity of the first author was supported by the Program (Vanvitelli per la Ricerca: VALERE) 2019 financed by the University of Campania ``L. Vanvitell''.}}
\date{}
\begin{document}
\markboth{\footnotesize{ F. Crispo and P.
Maremonti}} {\footnotesize{
A new proof of existence in the $L^3$-setting
}}
\maketitle
 \par\noindent
 \vskip -0.7true cm\noindent
\newcommand{\red}{\protect\bf}
\renewcommand\refname{\centerline
{\red {\normalsize \bf References}}}
\newtheorem{ass}
{\bf Assumption}[section]
\newtheorem{defi}
{\bf Definition}[section]
\newtheorem{tho}
{\bf Theorem}[section]
\newtheorem{rem}
{\sc Remark}[section]
\newtheorem{lemma}
{\bf Lemma}[section]
\newtheorem{coro}
{\bf Corollary}[section]
\newtheorem{prop}
{\bf Proposition}[section]
\renewcommand{\theequation}{\thesection.\arabic{equation}}
\setcounter{section}{0}
\numberwithin{equation}{section}
{\bf Abstract} - We investigate on the existence of solutions with initial datum $U_0$ in $L^3$. Our chief goal is to establish the existence interval $(0,T)$ uniquely considering the size and the absolute continuity of $|U_0(x)|^3$.  
\section{Introduction}\label{intro}
This note concerns the 3D-Navier-Stokes Cauchy problem:
 \be\label{NSC}\ba{l}U_t+U\cdot
\nabla U+\nabla\pi=\Delta
U,\;\nabla\cdot
U=0,\mbox{ in }(0,T)\times\R^3,\\
U=U_0\mbox{ on
}\{0\}\times\R^3.\ea\ee     
where $U$ and $\pi$ stand, respectively, for the unknown kinetic and pressure fields of an incompressible viscous fluid,
 $U_t:=
\frac\partial{\partial t}U$  and
 $U\cdot\nabla U:=
U_k\frac\partial{\partial x_k}U$. \par 
We look for a result of existence and uniqueness with an initial datum in $L^3(\R^3)$ and divergence free. It is known that this kind of result is not new. Indeed,  there is a wide literature on it, with a first contribution  due to T. Kato in \cite{Kt}.    Moreover, 
the $L^3$-metric of the existence class belongs to the set of scaling invariant metrics, this concept is meant as defined  in \cite{CKN}. \par This note does not aim  at giving an original result of existence in the $L^3$-class, but its interest is a little more specific, in a sense that we attempt to clarify below.
\par  As far as we know, considering   a scaling invariant $X$-norm for an   initial datum $U_0$,  the existence interval is global (in time) if $\dm U_0\dm_X$ is sufficiently small (in this regard we point out a recent contribution in weighted spaces with increasing weight \cite{GM}); otherwise, without any restriction, one proves the existence on some  interval  $(0,T)$, but no connection is given between the size of $\dm U_0\dm_X$ and a dimensionless size of $T$. Actually, the interval $(0,T)$ is  determined, by means of different strategies,  with the aid of other metrics and, as matter of fact, it is  deduced with respect to another metric. \par In this connection, the recent paper \cite{CMKC} seems to be an exception. It is employed the dimensionless weighted functional $\dm U_0\dm_{wt}^2:=\sup_{x}\intl{\R^3}\frac {U_0^2(y)}{|x-y|}dy$ and, in the set $L^2_{wt}$, where $\dm \cdot\dm_{wt}<\infty$,   the subset of the so called Kato class $K_3$ is considered. This special set of initial data furnishes, for the first time,  the interval of existence with a dimensionless size which involves properties of data, specifically the ones of the Kato class. We do not give further details of the results that arise, as they are a part of the wider set of results of this note, where we consider $\dm U_0\dm_{L^3}$  in place of $\dm U_0\dm_{wt}$. In this regard, we recall that  $\dm \cdot\dm_{wt}$ is not equivalent to the $L^3$-norm. We stress that we could also consider the $n$-dimensional Navier-Stokes Cauchy problem, of course considering the $L^n$-setting. 
\par In addition to the purpose of providing a new sufficient condition in $L^3$ (scaling invariant norm)  to establish the existence interval $(0,T)$, the results of this note are the starting point for a forthcoming paper, concerning the same questions but in the case of the initial boundary value problem in $(0,T)\times\O$, where $\po$ is assumed a sufficiently regular compact set, or $\O$ is the half-space.\par Here, we simply argue on the first question, and we only give a hint for the second question.\par In paper \cite{CMKC}, an element of the Kato class enjoys the following property:
$$\lim_{\rho\to0}\sup_x\intl{|x-y|<\rho}\frac{|U_0(y)|^2}{|x-y|}dy=0\,.$$ This limit property  allows us to state that $$\lim_{t\to0}t^\frac12\dm U^0(t)\dm_\infty=0\,,$$ where $U^0(t,x)$ is the solution to the heat equation with initial datum $U_0$.   Actually, setting  $\dm U_0\dm_{K^\rho}=\sup_x\intl{|x-y|<\rho}\frac{|U_0(y)|^2}{|x-y|}dy$, uniformly in $x$ we get
\be\label{LIC}\ba{ll}t^\frac12| U^0(t,x)|\hskip-0.2cm&\displ\leq \intl{|x-y|<\rho}H(x-y,t)|U_0(y)|dy+\intl{|x-y|>\rho}H(t,x-y)|U_0(y)|dy\VSE\leq c\dm U_0\dm_{K^\rho}+c\,\exp\big[-\sfrac{\rho^2}t\big]\dm U_0\dm_{wt}\,.\ea\ee This estimate is   the key tool to discuss the local existence of the solution to the  integral  equation given by means of the Oseen tensor  (see equation \rf{OL}), without requiring auxiliary conditions. Thanks to this estimate we are able to avoid  a time parameter to make coercive the integral part related to the  convective term. \par Here the strategy is the same. We replace the property of the Kato class with the absolute continuity of the integral, in particular of $|U_0(y)|^3\in L^1(\R^3)$. So that, \rf{LIC} is substituted by the following:
\be\label{LIC-LT}t^\frac12| U^0(t,x)|\leq  c\dm U_0\dm_{L^3(B(x,\rho))} +c\, \exp\big[-\frac{\rho^2}t\big]\dm U_0\dm_{3}\,,\ee where the first term  on the right-hand side,  for a suitable $\rho>0$, satisfies the absolute continuity property uniformly in $x\in\R^3$. However, in order to complete the relation which ensures the convergence of successive approximations (or the contraction principle) related to  the ``integral Oseen equation'', we need a more complete metric. For this task, for $\rho>0$, we set 
\be\label{TM}\dm u\dm_{3,\rho}^3:=\displ \sup_{\R^3}\!\intl{B(x,\rho)}\! |u(y)|^3\, dy, \ee
\be\label{TMM}
\displ\tm u\tm_{(t,\rho)}:=\sup_{(0,t)}\tau^\frac12\dm u(\tau)\dm_\infty+\sup_{(0,t)}\dm u(\tau)\dm_{3,\rho}+\frac {t^\frac12\hskip-0.15cm}{\rho}\hskip0.15cm\sup_{(0,t)}\dm u(\tau)\dm_{3}
\,.\ee 
 The second term in \rf{TMM}, roughly speaking, is preserved by the absolute continuity property of the initial datum, while the last term has the  dimensionless quantity  $\frac{t^\frac12}\rho$ as a weight, hence can be chosen suitably small. This strategy with weight for the norm in $L^3$ is practicable, as its estimate in the integral equation is linear in the $L^3$-norm, in particular in the convective term it is multiplied by the weighted (in time)  $L^\infty$-norm.   Since, for each norm in  \rf{TMM}, the convective term admits  an estimate in terms of $\tm U\tm_{(t,\rho)}^2\,,$
we conclude, with absolute  constants $c_0$ and $c_1$, with an estimate of the kind
$$\tm U\tm_{(t,\rho)}\leq c_0\Big[\dm U_0\dm_{3,\rho}+\frac{t^\frac12}\rho\dm U_0\dm_3\Big]+c_1\tm U\tm^2_{(t,\rho)}\,.$$
This last is the right one, as the term in square brackets can be chosen suitably to realize the wanted estimates, in particular the one of the existence. 
\par In this note we further look for some proprieties (see Corollary\,\ref{LQ}) of the solution  to the Cauchy problem \rf{NSC}. In some sense, these properties are close to the ones of the solution to the Stokes problem, and are crucial to obtain the same result in the case of the IBVP. Actually we are able to achieve the solution to the IBVP by employing the same smallness condition of the term $A(\rho,t)$ (see \rf{NIEQ}) that we require  for the existence of solutions to the Cauchy problem.
\vskip0.1cm\par We denote by $\mathscr C_0(\R^3)$ the subset of $C_0^\infty(\R^3)$
whose elements are  divergence free, and by $J^3(\R^3)$ the closure of $\mathscr C_0(\R^3)$ in $L^3(\R^3)$. 
Then we can state our chief results as follows:\begin{tho}\label{EKC}{There exist  an absolute constant $C>0$ such that  
 for all $U_0\in J^3(\R^3) $  there exists $T:=T(U_0)$, defined as 
 $$T(U_0)\!:=\displ\sup_{\rho>0}t(\rho), \mbox{with }  t(\rho):=\displ\sup\Big\{t>0:\  \Big[\dm U_0\dm_{3,\rho}+(e^{-\frac{\rho^2}{8t}} +\frac {t^\frac12\hskip-0.15cm}\rho\hskip0.15cm)\dm U_0\dm_{3}\Big]<\frac{1}{4C}\Big\}\,, $$
  such that problem \rf{NSC} has a solution $(U,\pi)$ on $(0,T)\times\R^3$ enjoying the properties \be\label{EKC-I}\ba{l}\mbox{for all } t\in(0,T(U_0)),\,\theta\in[0,1),\VS U\in  C^{2,\theta}(\R^3) \mbox{ and }U_t,\,D^2U\in C^{0,\frac\theta2}((\eta,T(U_0))\times\R^3), \VS
\sup_{(0,T)}t^\frac12\dm U(t)\dm_\infty\leq c\,\dm U_0\dm_{3},\ea\ee   
  with\be\label{EKC-L}\ba{c}\displ\lim_{t\to0}t^\frac12\dm U(t)\dm_\infty=0,	\quad \lim_{t\to0}t\dm \n U(t)\dm_\infty=0,\;,\VS \dm U(t)\dm_3\leq  c\dm U_0\dm_3\mbox{ for all }t\in[0,T(U_0))\,, \mbox{ and} \quad\lim_{t\to0}\|U(\tau)-U_0\|_3=0\,,\VS \displ t^\frac12\dm \pi(t)\dm_{3}\leq  c\,\dm U_0\dm_3\,,\;\pi\in C^{1,\theta}(\R^3)\mbox{ for all }t\in[0,T(U_0))\,.\ea \ee   Finally, if  the norm $\dm U_0\dm_{3}$ is suitably  small, then the above results hold for all $t>0$. }\end{tho} 
\begin{coro}\label{LQ}{Let $(U,\pi)$ the solution to problem \rf{NSC} stated in Theorem \ref{EKC}. Then,  for any $q\in (3,+\infty)$,  the following properties hold for $t\to 0$
 \be\label{SLQ}
 \dm U(t)\dm_q =o(t^{-\mu}),\ \dm \n U(t)\dm_q= o(t^{-\frac 12-\mu}),  \ 
\dm \n\n U(t)\dm_q= o(t^{-1-\mu}), \ee
with $
 \mu:=\frac {q-3}{2q}, $ and \be\label{SLIDS}\lim_{t\to0}t^{\frac32} \dm \n\n U(t)\dm_\infty=0\,.\ee}
\end{coro}
\begin{rem}{\rm 
We like to point out that we get a more detailed estimate than \rf{SLQ}. Actually,  
for any $q\in [3,+\infty)$ and for $ \mu:=\frac {q-3}{2q}$,  we have
  $$\ba{ll}\displ 
  t^{\mu}\dm U(t)\dm_q \leq c\dm U_0\dm_3^\frac3q A^\frac{q-3}q(\rho,t), \VS
  t^{\frac 12+\mu} \dm \n U(t)\dm_q \leq 
  c\dm U_0\dm_3^\frac3q \Big[A^\frac{q-3}q(\rho,t)+A^3(\rho,t)\Big], \VS
t^{1+\mu}\dm \n\n U(t)\dm_q\leq c\dm U_0\dm_3^\frac3q \Big[A^\frac{q-3}{q}(\rho,t)+A^6(\rho,t)\Big],\ea$$
}\end{rem}
for any $t\in [0,T)$, with $A(\rho,t)$ suitable function defined in \rf{NIEQ}. 
\begin{tho}[\bf Uniqueness]\label{PU}{\sl For all $U_0\in J^3(\R^3)$ a  solution to problem \rf{NSC} in the class of solutions $U\in L^\infty(0,T; J^3(\R^3))$,  satisfying \be\label{LPB}\ba{ll}\displaystyle  \lim_{t\to0}t^\frac12\dm U(t)\dm_\infty=0\,,\VS\displaystyle
\lim_{t\to 0}(U(t), \psi)=(U_0, \psi), \ \forall \psi\in \mathscr{C}_0(\R^3),\ea\ee is unique.}\end{tho}
This paper is organized as follows. In sect.2 we introduce some preliminary results.
In sect.3 we furnish estimates on the approximating sequence of solutions.  In sect.4 we prove the existence and uniqueness results (Theorem \ref{EKC} and Theorem \ref{PU}) and the $L^q$ limit properties of Corollary \ref{LQ}. 
\section{Preliminary results} We look for a solution to the integral equation
\be\label{OL}U(t,x)=H*U_0(t,x)-\n_xE*(U\otimes U)(t,x)\,,\,\mbox{ for all }(t,x)\in(0,T)\times\R^3\,,\ee
where $H(t,z):=(4\pi t)^{-\frac32}\exp[-{|z|^2}/{4t}]$ is the fundamental solution of the heat equation and $E(s,z)$ is the Oseen tensor, fundamental solution of the Stokes system, with components $$\ba{l}E_{ij}(s,z):=-H(s,z)\delta_{ij}+D_{z_iz_j}\phi(s,z)\,,\VS \phi(s,z):=\mathscr E(z)s^{-\frac32}\intll0{|z|} \exp[-{a^2}/{4s^\frac32}]da\,,\ea$$
where $\mathscr E$ is the fundamental solution of the Laplace equation. For the Oseen tensor the following estimates hold (see \cite{Kn}, estimates (VI) and (VIII) on pages 215 and 216, or  \cite{Os}):
\be\label{O-Med} \intl{\R^3}D_s^k D^h_z E(s,z)\,dz=0, \mbox{ for all }s>0\,,\ee
\be\label{O-I}|D_s^k D^h_z E(s,z)|\leq c(|z|+s^\frac12)^{-3-h-k},\mbox{ for all }s>0 \mbox{ and }z\in \R^3\,.\ee
For all $\theta\in(0,1)$, by the symbol $[g(t)]_{\theta}$ we denote the H\"older semi-norm. 
\par We also recall
\be\label{O-II}\ba{l}[D^h_zE(s)]_{\theta}\leq c \big[(|z|+s^\frac12)^{-(3+h+1)\theta}+(|\ov z|+\!s^\frac12)^{-(3+h+1)\theta}\big]\VS\hskip4.6cm\times\big[(|z|+s^\frac12)^{-(3+h)(1-\theta)}+(|\ov z|+\!s^\frac12)^{-(3+h)(1-\theta)}\big]\,,\VS [D^k_tE(z)]_{\frac\theta2}\leq c \big[(|z|+s^\frac12)^{-(3+h+1)\theta}+(|z|+\ov s^\frac12)^{-(3+{k}+1)\theta}\big]\VS\hskip4.6cm\times\big[(|z|+s^\frac12)^{-(3+k)(1-\theta)}+(|z|+\!\ov s^\frac12)^{-(3+h)(1-\theta)}\big]\,,\ea\ee where, for $h=\alpha_1+\alpha_2+\alpha_3$, $D^h_z$ denotes partial derivatives  with respect to $z_i$-variable $\alpha_i$ times, $i=1,2,3$. \par
\begin{lemma}\label{KRHOH}{\sl Let $\dm a\dm_{3,\rho}<\infty$. Then  for the   convolution product $H*a$ we get \be\label{KR-I}\dm H*a(t)\dm_{3,\rho}\leq  \dm a\dm_{3,\rho}\,,\mbox{ for all }t>0\,.\ee  }\end{lemma}
\P The result follows from a direct application of Minkowski's inequality. \chiu
\begin{lemma}\label{LI}{\sl Let $a\in L^{3}(\R^3)$. Then for the   convolution product $H*a$   we get 
\be\label{LI-I}\ba{c}t^\frac12\dm H*a(t)\dm_\infty+ t \dm \n H*a(t)\dm_\infty+t^\frac 32 \Big[\dm D_t H*a(t)\dm_\infty+\dm \nabla\nabla H*a(t)\dm_\infty\Big]
\VS \hfill \leq  h_0\dm a\dm_{3,\rho}+h_1\, e^{-\rho^{2}/8t}\dm a\dm_{3}\,,\mbox{ for all }\rho>0 \mbox{ and }t>0\,,
\ea\ee
with $h_0$ and $h_1$ positive constants. }\end{lemma}
 \P By the definition of the heat kernel and applying H\"older's inequality, 
 we get $$\ba{ll}| H*a(t,x)|\hskip-0.2cm&\displ\leq \!\!\intl{B(x,\rho)}H(t,x\!-\!y)|a(y)|\hskip0.15cmdy+\!\!\intl{|x-y|>\rho}H(t,x\!-\!y)|a(y)|\hskip0.15cmdy\VS&\displ\leq \Big[\!\!\intl{B(0,\rho)}\frac{e^{-\frac 32\frac{|z|^2}{4t}}}{(4\pi t)^{\frac 32}}dz\Big]^\frac23\dm a\dm_{3,\rho}+e^{-\frac {\rho^2}{8t}}\Big[\!\!\intl{|z|>\rho}  \frac{e^{-\frac 32\frac {|z|^2}{8t}}}{(4\pi t)^{\frac 32}}dz\Big]^\frac23\dm a\dm_{3}\VS&\displ \leq  t^{-\frac 12}\big[h_0\dm a\dm_{3,\rho}+h_1 e^{-\frac {\rho^2}{8t}}\dm a\dm_{3}\big]\,,
 \ea$$ where $h_0$ and $h_1$ are positive constants independent of $t$ and $\rho$. The other estimates follow by the same calculations, recalling that 
$$|D_s^k D^h_z H(s,z)|\leq c s^{\frac{-3-h}{2}-k} e^{-|z|^2/4t},\mbox{ for all }s>0 \mbox{ and }z\in \R^3\,. $$
 \chiu
\begin{lemma}\label{NLII}{\sl Let $\sup_{(0,T)}\!\!\big[\dm a(t)\dm_{3}+\dm b(t)\dm_{3}\big]\!<\infty$. Then there exists a constant $c$ independent of $a$ and $b$ such that, for $k>0$, \be\label{NLII-I}\hskip-0.1cm\ba{ll}\dy
\intll 0{\frac t2}\intl{\R^3}\frac{|a(\tau,y)||b(\tau,y)|}{(|x-y|+(t-\tau)^\frac 12)^{3+k}} \hskip0.1cm dy d\tau
\leq \displ ct^{-\frac{k}{2}}\Big[\sup _{(0,t)} \dm |a(\tau)||b(\tau)|\dm_{\frac32,\rho}\!\VS \hskip2cm\dy + t\rho^{-2}\!\sup _{(0,t)}\dm |a(\tau)||b(\tau)|\dm_{\frac32}\Big], \mbox{ for all }\rho>0\mbox{ for all  }t\in (0,T)\,.
\ea\ee}\end{lemma}
\P By H\o lder's inequality and the hypotheses,  $$\ba{ll}\dy \intll 0{\frac t2}\intl{\R^3}\frac{|a(\tau,y)||b(\tau,y)|}{(|x-y|+(t-\tau)^\frac 12)^{3+k}} \hskip0.1cm dy d\tau
\VS \leq c\intll0{\frac t2} \intl{B(x,\rho)}\frac{|a(\tau,y)||b(\tau,y)|}{(|x-y|+(t-\tau)^\frac 12)^{3+k}}dyd\tau+\intll0{\frac t2} \intl{|x-y|>\rho}\frac{|a(\tau,y)||b(\tau,y)|}{(|x-y|+(t-\tau)^\frac 12)^{3+k}}dyd\tau\VS\displ\leq c\intll0{\frac t2}(t-\tau)^{-1-\frac{k}{2}}\dm |a(\tau)||b(\tau)|\dm_{\frac32,\rho}d\tau+
c\rho^{-2}\intll0{\frac t2}(t-\tau)^{-\frac{k}{2}}\dm |a(\tau)||b(\tau)|\dm_{\frac32}d\tau\VS\leq ct^{-\frac{k}{2}}\sup _{(0,t)}\dm |a(\tau)||b(\tau)|\dm_{\frac32,\rho}+ct^{1-\frac{k}{2}}\rho^{-2}\sup _{(0,t)}\dm |a(\tau)||b(\tau)|\dm_{\frac32},\ea$$ that is  \rf{NLII-I}.\chiu
\begin{lemma}\label{LII}{\sl Let $\sup_{(0,T)}\!\!\big[t^\frac12\dm a(t)\dm_\infty+t^\frac12\dm b(t)\dm_\infty\big]\!<\infty$ and $\sup_{(0,T)}\!\!\big[\dm a(t)\dm_{3}+\dm b(t)\dm_{3}\big]\!<\infty$. Then there exists a constant $c$ independent of $a$ and $b$ such that
\be\label{LII-I}\hskip-0.1cm\ba{l}t^\frac12\dm \n E*\!(a\otimes b)(t)\dm_\infty \!\leq \displ c\Big[\sup _{(0,t)}\tau\dm |a(\tau)|| b(\tau)|\dm_\infty\!\VS\hskip3.5cm+\sup _{(0,t)} \dm |a(\tau)|^\frac12|b(\tau)|^\frac12\dm_{3,\rho}^2\!+ t\rho^{-2}\!\sup _{(0,t)}\dm |a(\tau)|^\frac12|b(\tau)|^\frac12\dm_{3}^2\Big], \VS\hskip 7cm\mbox{ for all }\rho>0\mbox{ for all  }t\in (0,T)\,.\ea\ee}\end{lemma}
\P  Via formulae \rf{O-I}  we get
$$\ba{ll}\displ|\n E*(a\otimes b)(t,x)|\hskip-0.2cm&\displ\leq  \intll 0{\frac t2}\intl{\R^3}\frac{|a(\tau,y)||b(\tau,y)|}{(|x-y|+(t-\tau)^\frac 12)^4\hskip-0.1cm}\hskip0.1cmdyd\tau+
\intll {\frac t2}t\intl{\R^3}\frac{ |a(\tau,y)||b(\tau,y)|}{(|x-y|+(t-\tau)^\frac 12)^{4}\hskip-0.1cm}\hskip0.1cmdyd\tau\displ\VSE=:I_1(t)+I_2(t)\,.\ea$$
By our hypotheses we get
$$I_2(t)\leq c\intll {\frac t2}t\frac 1\tau\sup \tau\dm |a(\tau)||b(\tau)|\dm_\infty\intl{\R^3}(|z|^2+t-\tau)^{-2}dzd\tau\leq ct^{-\frac12}\sup _{(\frac t2,t)}\tau\dm |a(\tau)||b(\tau)|\dm_\infty\,.$$
For $I_1$ we use estimate \rf{NLII-I} with $k=1$, and we get $$I_1(t)\leq ct^{-\frac12}\sup _{(0,t)}\dm |a(\tau)|^\frac12|b(\tau)|^\frac12\dm_{3,\rho}^2+ct^\frac12\rho^{-2}\sup _{(0,t)}\dm |a(\tau)|^\frac12|b(\tau)|^\frac12\dm_{3}^2.$$ From the previous we arrive at \rf{LII-I}.\chiu
\begin{lemma}\label{KRHOE}{\sl Let $\sup_{(0,T)}\big[ t^\frac12\dm a(t)\dm_\infty+\dm b(t)\dm_{3,\rho}\big]<\infty$. Then there exists a constant $c$ independent of $a(t,x)$ and $b(t,x)$  such that \be\label{KR-II}\dm \n E*(a\otimes b)(t)\dm_{3,\rho}\leq c\sup _{(0,t)}\tau^\frac12\dm a(\tau)\dm_\infty \dm b(\tau)\dm_{3,\rho}, \mbox{ for all }\rho>0\mbox{ and }t\in (0,T)\,.\ee}\end{lemma}
\P Set $\xi=y-z$ in the convolution product. We have   $$\dm \n E*(a\otimes b)(t)\dm_{3,\rho} =\Big[\intl{B(x,\rho)}\Big[\intll0t\intl{\R^3}\n E
(t-\tau,\xi)\cdot (a(\tau,y-\xi)\otimes b(\tau,y-\xi))d\xi d\tau\Big]^3dy\Big]^\frac 13$$
Employing Minkowski's inequality, then our hypotheses and estimate \rf{O-I} for the Oseen tensor, we find 
$$\ba{ll}\dy\vspace{1ex}
\dm \n E*(a\otimes b)(t)\dm_{3,\rho}&\dy\leq \intll0t\intl{\R^3}|\n E(t-\tau,\xi)|\Big[\intl{B(x,\rho)}|a(\tau,y-\xi)|^3|b(\tau,y-\xi)|^3 dy\Big]^\frac13d\xi d\tau\\ &\dy \leq c\sup _{(0,t)}\tau^\frac12\dm a(\tau)\dm_\infty\dm  b(\tau)\dm_{3,\rho}\intll0t(t-\tau)^{-\frac12}\tau^{-\frac12}d\tau\,.\ea$$ that  gives \rf{KR-II}. 
\chiu
\begin{lemma}\label{KWHOLE}{\sl Let $\sup_{(0,T)}\big[t^\frac12\dm a(t)\dm_\infty+\dm b(t)\dm_{3}\big]<\infty$. Then there exists a constant $c$ independent of $a(t,x)$ and $b(t,x)$  such that \be\label{KW-I}\dm \n E*(a\otimes b)(t)\dm_{3}\leq c\sup _{(0,t)}\tau^\frac12\dm a(\tau)\dm_\infty \dm b(\tau)\dm_3, \mbox{ for all }t\in (0,T)\,.\ee}\end{lemma}
\P The proof is analogous to the one of the  previous lemma. Hence it is omitted.\chiu
\begin{lemma}\label{HHST}{\sl In the hypotheses of Lemma\,\ref{LI} and Lemma\,\ref{LII} the convolution products  $\n H*a$ and $\n E*(a\otimes b)$ are H\"oder continuous functions in $x\in \R^3$, with  exponent $\theta\in[0,1)$.
In particular, we get \be\label{HV} \hskip-0.2cm\ba{ll}\mbox{\hskip0.8cm\large$\frac{|H*a(t,x)-H*a( t,\ov x)|}{|x-\ov x|^\theta} $}&\hskip-1cm\leq c\hskip0.05cmt^{-\frac{1+\theta}2}
\big[h_0\dm a\dm_{3,\rho}+h_1 e^{-\frac {\rho^2}{8t}}\dm a\dm_{3}\big]\,,
\vspace{4pt}\\
\mbox{\hskip0.8cm\large$\frac{|\n H*a(t,x)-\n H*a( t,\ov x)|}{|x-\ov x|^\theta} $}&\hskip-0.4cm\leq c\hskip0.05cmt^{-1-\frac{\theta}2}\big[h_0\dm a\dm_{3,\rho}+h_1 e^{-\frac {\rho^2}{8t}}\dm a\dm_{3}\big]\,,\vspace{4pt}
\\\mbox{\large$\frac{|\n\hskip-0.02cm E*\hskip-0.02cm(a\otimes b)(t,x)-\!\n\hskip-0.02cm E*\hskip-0.02cm(a\otimes b)( t,\ov x)|}{|x-\ov x|^\theta}$}&\hskip-0.2cm\leq c\hskip0.05cmt^{-\frac{1+\theta}2}\big[\mbox{$\displ\sup _{(0,t)}$}\hskip0.03cm\tau\dm a(\tau)b(\tau)\dm_\infty\!+\mbox{$\displ\sup _{(0,t)}$}\dm |a(\tau)b(\tau)|^\frac12\dm_{3, \rho}^2\VS&\hskip-0.5cm \dy\hskip2cm +
\frac{t}{\rho^2} \dm |a(\tau)b(\tau)|^\frac12\dm_{3}^2
\big]\hskip-0.025cm,\mbox{ for all } \rho>0,\ea\ee
 with $h_0$ and $h_1$ positive constants independent of $t$ and $\rho$.  }\end{lemma}
\P 
 The first two estimates follow  applying the Lagrange theorem and employing the $L^\infty$ estimates of Lemma\,\ref{LI} for the convolution products  $ H*a$, $\n H*a$ and $\n \n H* a$. Hence we limit ourselves to prove estimate \rf{HV}$_{3}$. 
 From properties \rf{O-II} for the Oseen tensor $E$,  we get
$$\ba{ll}| \n E*(a\otimes b)(t,x)-\n E*(a\otimes b)( t,\ov x)|\!\VS \displ \leq c|x-\ov x|^\theta\!\intll0{\frac t2}\!\intl{\R^3}|a(\tau,y)||b(\tau,y)|\Big[\frac1{(|x-y|+(t-\tau)^\frac12)^{4+\theta}}+\frac 1{(|\ov x-y|+(t-\tau)^\frac12)^{4+\theta}}\Big]dyd\tau\VS +
\displ c|x-\ov x|^\theta\!\intll{\frac t2}t\! \intl{\R^3}|a(\tau,y)||b(\tau,y)|\Big[\frac1{(|x-y|+(t-\tau)^\frac12)^{4+\theta}}+\frac 1{(|\ov x-y|+(t-\tau)^\frac12)^{4+\theta}}\Big]dyd\tau\VS
=:I_1+I_2\ea$$
For $I_2$ we get
$$\ba{ll}I_2\leq \displ c|x-\ov x|^\theta\sup _{(0,t)}\hskip0.03cm\tau\dm a(\tau)b(\tau)\dm_\infty
\intll{\frac t2}t\tau^{-1}\!\!\intl{\R^3}\Big[\frac1{(|x-y|+(t-\tau)^\frac12)^{4+\theta}}\VS \hfill\displ+\frac 1{(|\ov x-y|+(t-\tau)^\frac12)^{4+\theta}}\Big]dyd\tau
\VS
\displ \hskip0.5cm \leq c |x-\ov x|^\theta t^{-\frac12-\frac\theta 2}\sup _{(0,t)}\hskip0.03cm\tau\dm a(\tau)b(\tau)\dm_\infty.\ea.$$
For $I_1$ we apply estimate \rf{NLII-I} with exponent $k=1+\theta$ 
and we easily obtain
$$\ba{ll}I_1\displ\leq  c|x-\ov x|^\theta t^{-\frac 12-\frac \theta 2}\Big[\sup _{(0,t)}\dm |a(\tau)|^\frac 12	|b(\tau)|^\frac12\dm_{3,\rho}^2+
 t \rho^{-2}\sup _{(0,t)}\dm |a(\tau)|^\frac 12	|b(\tau)|^\frac12\dm_{3}^2\Big].
\ea$$
From the previous estimates, for all  $t\in(0,T)\times\R^3$ we easily get \rf{HV}$_2$.  \chiu
\begin{lemma}\label{HHTT}{\sl In the hypotheses of Lemma\,\ref{LI} and Lemma\,\ref{LII} the convolution products  $\n H*a$ and $\n E*(a\otimes b)$ are H\"older continuous functions in $t\in (0,T)$, with  exponent $\theta\in[0,\frac12)$. }
\end{lemma}
\P The proof could be obtained arguing as for the H\"older property with respect to the space variable. On the other hand, for our aims we don't need estimates of the kind \rf{HV}. Hence we omit further details.
\chiu
 We set $e_i(t,z):=\n E_{i}(t,z)$, $i\in \{1,2,3\}$, and, for some tensor field $w$, we set $$W(t,x):=\intll0t\intl{\R^3}e_i(t-\tau,x-y)\cdot w(\tau,y)dyd\tau\,.$$ \begin{lemma}\label{KG}{\sl Let $w(t,x)\in L^\infty(0,T; L^\frac 32(\R^3))$ and   $w(t,x)\in C^{0,\theta}(\R^3)$ for all $t\in (0,T)$  with $$\sup_{(\frac t2,t)}\tau^{1+\frac\theta2}[w(t)]_\theta<\infty,\mbox{ for all }t\in(0,T)\,,$$ then, for all $\ov\theta<\theta$,  we get 
\be\label{GW}\ba{l}t\dm \n W(t)\dm_\infty+t^{1+\frac{\ov\theta}2}[\n W(t)]_{\ov\theta}\leq c\displ\sup_{(\frac t2,t)}\tau^{1+\frac\theta2}[ w(\tau)]_\theta\VS\hskip6cm+ c\sup_{(0,\frac t2)}\hskip-0.05cm\Big[\dm  w(\tau)\dm_{\frac32,\rho}\hskip-0.05cm+\frac t{\rho^2} \dm w(\tau)\dm_\frac32\Big],\ea \ee \mbox{for all }$t\in (0,T)$\,. Moreover, if $\n\cdot w\in C^{0,\theta}(\R^3)$ for all $t\in (0,T)$,  with $$\sup_{(0,\frac t2)}\tau^{\frac32+\frac\theta2}[\n\cdot w(\tau)]_\theta<\infty\,,\mbox{ for all }t\in(0,T)\,,$$ then, for all $\ov\theta<\theta$,  we get 
\be\label{GGW}\ba{l}t^\frac32\dm \n\n W(t)\dm_\infty\!+t^{\frac32+\frac{\ov\theta}2}[\n\n W(t)]_{\ov\theta}\leq c\displ\sup_{(\frac t2,t)}\tau^{\frac32+\frac\theta2}[\n\hskip-0.05cm\cdot\hskip-0.05cm w(\tau)]_\theta\VS\hskip6.6cm+ c\hskip-0.05cm\sup_{(0,\frac t2)}\Big[ \dm  w(\tau)\dm_{\frac32,\rho}\hskip-0.05cm+\frac t{\rho^2} \dm  w(\tau)\dm_\frac32\Big],\ea\ee \mbox{for all }$t\in (0,T)$\,.}\end{lemma} \P We set $$W_\vep(t,x) := \int_0^{t-\vep}\!
\intl{\R^3} e_i(t-\tau, x-y)\cdot w(\tau,y)dyd\tau.$$
By using the H\o lder property of $w$, a classical argument, ensures the existence of
 $$
 \lim_{\vep\to 0}\n W_{\vep}(t,x)= \n W(t,x)$$
with \be\label{VNW}\n W(t,x):= \lim_{\vep\to0}\intl0^{t-\vep}
\intl{\R^3} \n e_i(t-\tau, x-y)\cdot w(\tau,y)dyd\tau. \ee
 Let us write $\n W_\vep(t,x)$ as follows 
$$\ba{ll}\dy \n W_\vep(t,x) = \intll0{\frac t2}\!
\intl{\R^3} \n_x e_i(t-\tau, x-y)\cdot w(\tau,y)dyd\tau\\\hskip3cm\displ +\intll{\frac t2}{t-\vep}\!
\intl{\R^3} \n_x e_i(t-\tau, x-y)\cdot(w(\tau,y)-w(\tau, x))\,dyd\tau=:I_1+I_2.\ea$$
By using property \rf{O-I} and Lemma\,\ref{NLII}, with $k=2$,  we find
$$\ba{ll} |I_1|\leq c\dy \intll0{\frac t2}\!\!
\intl{\R^3} \frac{|w(\tau, y)|}
{(|x-y|+(t-\tau)^{\frac12})^5}dy d\tau
\dy\leq   \frac ct\sup_{(0,\frac t2)}\Big[\dm w\dm_{\frac 32,\rho}+\frac t{\rho^2} \dm w(\tau)\dm_\frac 32\Big]\,,\mbox{ for all }t\in (0,T)\,.
\ea$$
By using property \rf{O-I} and the H\o lder property of $w$, for $I_2$ we get
$$\ba{ll} |I_2|&\displ\leq c\sup_{(\frac t2,t)}\tau^{1+\frac\theta2}[w(\tau)]_\theta\intll{\frac t2}{t-\vep}\!\!\!
\tau^{-1-\frac {\theta}{2}}\intl{\R^3} \frac{1}
{(|x-y|+(t-\tau)^{\frac12})^{5-\theta}}dy d\tau\VS &\leq 
\dy \frac {c}{t^{1+\frac \theta 2}}
\sup_{(\frac t2,t)}\tau^{1+\frac\theta2}[w(\tau)]_\theta\intll{\frac t2}{t-\vep}\!\!(t-\tau)^{-1+\frac \theta 2} d\tau\VS &\dy\leq \frac{c}{t}\sup_{(\frac t2,t)}\tau^{1+\frac\theta2}[w(\tau)]_\theta\,,\mbox{ for all }t\in (0,T).
\ea$$
Hence, uniformly in $\vep>0$, we arrive at 
$$t\dm \n W_\vep (t)\dm_\infty\leq  c\sup_{(0,\frac t2)}\tau^{1+\frac\theta2}[w(\tau)]_\theta +c\sup_{(\frac t2,t)}\Big[\dm w\dm_{\frac 32,\rho}+\frac t{\rho^2} \dm w(\tau)\dm_\frac 32\Big]\,,\mbox{ for all }t\in (0,T),$$ which leads to the $L^\infty$-estimate  enclosed in \rf{GW}.
\par
Now we show the H\"older property of $\n W$. 
 We set
$$\ba{l}| \n W_\vep(t,x)-\n W_\vep(t,\ov x)|\\\hskip2cm 
=\dy \bigg|\intll0{\frac t2}\!
\intl{\R^3}[ \n e_i(t-\tau, x-y)-\n e_i(t-\tau, \ov x-y)]\cdot  w(\tau,y)dyd\tau\\\hskip4.2cm +\displ\intll{\frac t2}{t-\vep}\!
\intl{\R^3} [\n e_i(t-\tau, x-y)- \n e_i(t-\tau, \ov x-\tau)]\cdot w(\tau, y)\,dyd\tau.
 \ea$$
From properties \rf{O-II} for the Oseen tensor $E$ and H\"older's assumption on $w$, for all $\theta'<\theta$,  we easily get
$$\ba{l}|I_1(t)| \displ \leq c|x-\ov x|^{\theta'}\int_0^\frac t2\!\!\intl{\R^3}|w(\tau,y)| \Big[\ssfrac1{(|x-y|+t^{\frac12})^{5+\theta'}}+\ssfrac 1{(|\ov x-y|+t^{\frac12})^{5+\theta'}}\Big]dyd\tau\,,\\|I_2(t)|\displ\leq \frac{c}{t^{1+\frac {\theta}{2}}} \displ\sup_{(\frac t2,t)}\tau^{1+\frac\theta2}[w(\tau)]_\theta \displ |x-\ov x|^{\theta'}\!\!\int_{\frac t2}^t \intl{\R^3}\Big[\ssfrac1{(|x-y|+t^{\frac12})^{5+\theta'-\theta}}+\ssfrac 1{(|\ov x-y|+t^{\frac12})^{5+\theta'-\theta}}\Big]dyd\tau	\,.\ea$$
By using Lemma \ref{NLII} with $k=2+\theta'$, for $I_1$ we find
$$|I_1|\leq 
c|x-\ov x|^{\theta'}
t^{-1-\frac{\theta'\hskip-0.1cm}{2}}\hskip0.1cm\Big[\sup _{(0,\frac t2)} \dm w\dm_{\frac32,\rho}+ t\rho^{-2}\!\sup _{(0,\frac t2)}\dm w(\tau)\dm_{\frac32}\Big]\,.
$$
For $I_2$ an integration furnishes
$$\ba{ll} \dy |I_2| \leq
\displ c|x-\ov x|^{\theta'} 
t^{-1-\frac{\theta'\hskip-0.1cm}{2}}\hskip0.1cm\sup_{(\frac t2,t)}\tau^{1+\frac\theta2}[w(t)]_\theta\,.
\ea$$
Since the estimates are uniform in $\vep>0$,  from the above estimates for $I_1$ and $I_2$ we obtain the thesis for all $\ov\theta<\theta'$. In order to prove \rf{GGW}, recalling the definition of $e_i$, then it is enough to consider that an integration by parts furnishes 
$$\ba{l} |\n\n W_\vep(t,x)|
=\dy \Big|\intll0{\frac t2}\!
\intl{\R^3} \n\n e_i(t-\tau, x-y)\cdot w(\tau,y)dyd\tau\\\hskip3.5cm -\displ\intll{\frac t2}{t-\vep}\!
\intl{\R^3} \n\n E_i(t-\tau, x-y)\cdot\n\cdot w(\tau, y)\,dyd\tau\Big|\,.
 \ea$$ After which, one employs the same arguments considered for $\n W$.
\chiu \begin{lemma}\label{GLT-L}{\sl If $t^\frac12 w,\,t\n\cdot w(t)\in L^\infty(0,T;L^3(\R^3))$, then, we get
\be\label{GLT}t^\frac12\dm \n W(t)\dm_3\leq c\Big[\sup_{(0,\frac t2)}\tau^\frac12\dm w(\tau)\dm_3+\displ\sup_{(\frac t2,t)}\tau\dm \n\cdot w(\tau)\dm_3\Big], \mbox{ for all }t\in (0,T)\,,\ee and if we also assume $t^\frac32\n\n\cdot w\in L^\infty(0,T,L^3(\R^3))$, then, we get \be\label{GGLT} t\dm \n\n W(t)\dm_3\leq c\Big[\sup_{(0,\frac t2)}\tau^\frac12\dm w(\tau)\dm_3+\sup_{(\frac t2,t)}\tau^\frac32\dm \n\n\cdot w(\tau)\dm_3,\mbox{ for all }t\in (0,T)\Big]\,.\ee}\end{lemma}\P We prove \rf{GLT}. Recalling the definition of $e_i$,
via an integration by parts, we have
\be\label{GLT-I}\ba{ll}\dy \dm \n W(t)\dm_3
=
\Big[\intl{\R^3}\bigg|\intll0{\frac t2}\intl{\R^3} \n e_i(y-z, t-\tau)\cdot w(z,\tau)dzd\tau dyd\tau
\\\displ\hskip4.3cm+ \intll{\frac t2}t\intl{\R^3} \n E_i(y-z, t-\tau)\cdot \n\cdot w(z,\tau)dzd\tau\bigg|^3dy\Big]^\frac 13.
\ea\ee
Applying for both terms Minkowski's inequality, estimate \rf{O-I} for the Oseen tensor, we find
$$\ba{ll}\dy \dm \n W(t)\dm_3\leq c \sup_{(0,\frac t2)}\tau^\frac 12\dm w(\tau)\dm_3\intll0{\frac t2}\tau^{-\frac 12} \intl{\R^3}|\n e_i(t-\tau, \xi)|d\xi d\tau\\\displ \hskip5cm+c \sup_{(\frac t2,t)}\tau\dm \n \cdot w(\tau)\dm_3\intll{\frac t2}t\tau^{-1}\intl{\R^3}| e_i(t-\tau, \xi)|d\xi d\tau 
\vspace{2pt}\\\displ\hskip1.6cm \leq ct^{-\frac 12} \Big[\sup_{(0,\frac t2)}\tau^\frac 12\dm w(\tau)\dm_3+ \sup_{(\frac t2,t)}\tau\dm \n \cdot w(\tau)\dm_{3}\Big], 
\ea$$
which proves \rf{GLT}. To prove \rf{GGLT}, it is enough to employ a further integration by  parts  on $(\frac t2,t)\times\R^3$, after which the argument lines are the same, so that we consider achieved the lemma.  
 \chiu
Let us consider the equation \be\label{Pres}(\pi,\Delta g)=-(a\otimes  u,\n\n g)\mbox{ for all }g\in C_0^\infty(\R^3).\ee   From 
the theory of singular integrals \cite{St}, one immediately gets the following result.  
\begin{lemma}\label{LPE-L}{\sl Let $a$ and $u$ be divergence free.  If $a\in L^\infty(\R^3)$ and $u\in L^3(\R^3)$, there exist constants $c$ independent of $a$ and $u$ such that for a solution $\pi$ to problem \rf{Pres}  the following holds: 
\be\label{LP-I} \dm \pi\dm_3\leq c\dm a\dm_\infty\dm u\dm_3\,.\ee
Further, if $a,u\in C^{1,\theta}(\R^3)$,  for some $\theta\in(0,1)$,  then $\pi\in C^{1,\theta}(\R^3)$.}\end{lemma}
\section{Properties of the approximating sequence of solutions}
We study the integral relation  \be\label{IEQ}U^m(t,x)=H*U_0(t,x)-\n_x E*(U^{m-1}\otimes U^{m-1})(t,x)\,.\ee 
\begin{lemma}\label{LIEQ}{\sl Let $U_0(x)\in L^3(\R^3)$. Set $U^0(t,x):=H*U_0$. Then there exist  constants $h_0$, $h_1$, $c_1$,  independent of $U_0$ and $m\in\N$, such that for the sequence \rf{IEQ} we get
\be\label{IEQ-I}\ba{l}\tm U^m\tm_{(t,\rho)}\leq (h_0+1)\dm U_0\dm_{3,\rho}+(h_1\, e^{-\rho^{2}/8t}+\frac {t^\frac12\hskip-0.15cm}{\rho}\hskip0.15cm)\dm U_0\dm_3+c_1\tm U^{m-1}\tm^2_{(t,\rho)}\,,\VS\hskip8cm \mbox{ for all }t>0\mbox{ and }\rho>0\,.\ea\ee}\end{lemma}
\P From definition \rf{IEQ}, by virtue of Lemma\,\ref{LI}-Lemma\,\ref{KWHOLE}, for all $t>0$ and $\rho>0$,  we get 
$$\ba{ll}\hskip0.1cms^\frac12\dm U^1(s)\dm_\infty\hskip-0.2cm&\leq \!
h_0\dm U_0\dm_{3,\rho} +h_1\, e^{-\rho^{2}/8s}\dm U_0\dm_{3}\VSE\hskip2cm+c\Big[\sup_{(0,s)}\tau^\frac12\dm U^0(\tau)\dm_\infty\hskip-0.05cm+\sup_{(0,s)} \dm U^0(\tau)\dm_{3,\rho}\!+ \mbox{\Large$\frac {s^\frac12\hskip-0.1cm}{\rho}$}\sup_{(0,s)}\dm U^0(\tau)\dm_{3}\Big]^2\!,\VS\hskip0.4cm\dm U^1(s)\dm_{3,\rho}\hskip-0.2cm&\leq \dm U_0\dm_{3,\rho}+\displ c\sup_{(0,s)}\tau^\frac12\dm U^0(\tau)\dm_\infty\dm U^0(\tau)\dm_{3,\rho}\VSE\leq \dm U_0\dm_{3,\rho}+\displ c\Big[\sup_{(0,s)}\tau^\frac12\dm U^0(\tau)\dm_\infty\hskip-0.05cm+\sup_{(0,s)} \dm U^0(\tau)\dm_{3,\rho}\!+ \mbox{\Large$\frac {s^\frac12\hskip-0.1cm}{\rho}$}\sup_{(0,s)}\dm U^0(\tau)\dm_{3}\Big]^2\!,\VS\dy \hskip0.4cm\dm U^1(s)\dm_{3}\hskip-0.2cm&\dy\leq \dm U_0\dm_{3}+ 
c\sup_{(0,s)}\tau^\frac12\dm U^0(\tau)\dm_\infty \dm U^0(\tau)\dm_3\,,\ea$$ where $c$ is a constant independent of $t,\,\rho$. Multiplying the last estimate for \mbox{\large$\frac {s^\frac12}{\rho}$}, and then increasing, we get
$$\mbox{\Large$\frac {s^\frac12}{\rho}$}\dm U^1(s)\dm_{3}\leq \mbox{\Large$\frac {s^\frac12}{\rho}$}\dm U_0\dm_{3}+ 
c\Big[\sup_{(0,s)}\tau^\frac12\dm U^0(\tau)\dm_\infty\hskip-0.05cm+\sup_{(0,s)} \dm U^0(\tau)\dm_{3,\rho}\!+ \mbox{\Large$\frac {s^\frac12\hskip-0.1cm}{\rho}$}\sup_{(0,s)}\dm U^0(\tau)\dm_{3}\Big]^2\,.$$Taking $\sup_{(0,t)}$    of the previous trilogy, then summing the first two  with the last one, recalling the definition of the functional $\tm\hskip0.05cm\cdot\hskip0.05cm\tm_{(t,\rho)}$, we arrive at
$$\ba{ll}\tm U^1\tm_{(t,\rho)}\leq \!(h_0+1)\dm U_0\dm_{3,\rho}+( h_1\,e^{-\rho^{2}/8t}+\mbox{\Large$\frac {t^\frac12}{\rho}$})\dm U_0\dm_{3}+3c\,\tm U^0\tm^2_{(t,\rho)}\,,\VS\hskip8cm \mbox{ for all }\rho>0\mbox{ and }t>0,\ea$$ with a constant $c$ independent of the datum $U_0$. So that, for $m=1$, \rf{IEQ} is well defined and estimate \rf{IEQ-I} is true. Then by induction   one proves the estimate for all $m\in\N$.
\chiu
We use the method of successive approximations. We show that the previous lemmas ensure boundedness and convergence of the approximating sequence of velocity fields $\{U^m\}$. Firstly we recall the following result.
\begin{lemma}\label{Sl}{\sl Let $\xi_0>0$ and $c>0$. Let $\{\xi_m\}$ be a non negative sequence of real numbers such that $$ \xi_m\leq \xi_0+c\xi_{m-1}^2\,.$$ Assume $1-4c\xi_0>0$ and $\xi_0\leq \xi$, where $\xi$ is the minimum solution of the algebraic equation $c\xi^2-\xi+\xi_0=0$. Then $\xi_{m-1}\leq \xi$ for all $m\in\N$.}\end{lemma}\P For the proof we refer to \cite{Sol77}. \chiu
\begin{lemma}\label{CONV}{\sl Let $\{U^m\}$ be the sequence  defined in \rf{IEQ} corresponding to $U_0\in J^3(\R^3)$. Then, there exists a $T(U_0)>0$ such that, for all $\eta$, the sequence strongly converges  in $C((\eta,T(U_0))\times\R^3),$ to a solution $U$ to \rf{OL}, and, for all $t\in(0,T(U_0))$, the sequence converges   to $U$ in $ J^3(\R^3)$. In particular we get, for a suitable $\rho$ and for all $t\in[0,T(U_0))$,
\be\label{CRE}\tm U\tm_{(t,\rho)}\leq \frac{2\big[(h_0+1)\dm U_0\dm_{3,\rho}+( h_1e^{-\rho^{2}/8t} +\frac {t^\frac12\hskip-0.15cm}{\rho}\hskip0.15cm)\dm U_0\dm_{3}\big]}{1+\big(1-4c_1\big[(h_0+1)\dm U_0\dm_{3,\rho}+(h_1e^{-\rho^{2}/8t}+\frac {t^\frac12\hskip-0.15cm}{\rho}\hskip0.15cm)\dm U_0\dm_{3}\big]\big)^\frac12}, \ee and
\be\label{LP1}\dm U(t)\dm_3\leq c \dm U_0\dm_3.\ee
Further \be\label{LP}\lim_{t\to0}\sup_{(0,t)}\tau^\frac12\dm U(\tau)\dm_\infty=0\,.\ee
}\end{lemma}
\P Since $U_0\in L^3(\R^3)$, for any $\vep\in (0,\frac{1}{4c_1(h_0+1)})$, there exists $\rho=\rho(U_0, \vep)$ such that 
$\dm U_0\dm_{3,\rho}<\vep$. For any such $\rho$, we denote by $t(\rho)$ the supremum of $t>0$ for which the following inequality holds
\be\label{EIEQ}  1-4c_1\big[(h_0+1)\dm U_0\dm_{3,\rho}+(h_1e^{-\frac{\rho^2}{8t}} +\frac {t^\frac12\hskip-0.15cm}\rho\hskip0.15cm)\dm U_0\dm_{3}\big]>0\,.
\ee
We observe that the definition of $t(\rho)$ is well posed, taking into account that, for any fixed $\rho>0$, the function in round brackets is a monotonic increasing function of $t$ that tends to zero as $t\to 0$. Finally we denote by $T(U_0)$  the supremum of   $t(\rho)$  for which \rf{EIEQ} holds. Then, by virtue of estimate \rf{IEQ-I} and applying Lemma\,\ref{Sl}, for a fixed $\rho$ and for any $t\in[0,T(U_0))$ and uniformly in $m\in\N$ we get\be\label{NIEQ}\tm U^m\tm_{(t,\rho)}\leq\!\frac{2\big[(h_0+1)\dm U_0\dm_{3,\rho}+(h_1e^{-\frac{\rho^2}{8t}} +\frac {t^\frac12\hskip-0.15cm}{\rho}\hskip0.15cm)\dm U_0\dm_{3}\big]}{1+\big(1-4c_1\big[(h_0+1)\dm U_0\dm_{3,\rho}+(h_1e^{-\frac{\rho^2}{8t}}  +\frac {t^\frac12\hskip-0.15cm}{\rho\hskip0.15cm})\dm U_0\dm_{3}\big]\big)^\frac12}\!\!=:A(\rho,t).\ee   
Estimate \rf{NIEQ} ensures that, for all $t\in[0,T(U_0)),$ the sequence $\{\tm U^m\tm_{(t,\rho)}\}$ is bounded.  
\par On the other hand, the validity of estimate \rf{IEQ-I}, for any $\rho>0$ and for any $t>0$, ensures
that the following property holds true:\begin{itemize}\item[\texttt{P:}]{\sl 
 for any  sequence $\{t_p\}\to0$,  one can construct a sequence $ \{\rho_p\}\to0$  such that  \rf{EIEQ} holds, and along these sequences we get 
$\lim_{p\to \infty} A(\rho_p, t_p)=0$ too\footnote{It is sufficient to choose the sequence $\{\rho_p\}$ such that $\rho_p\to 0 $ and $\frac{t_p^\frac 12}{\rho_p}=o(1)$. }. }\end{itemize}Therefore, we can again apply Lemma\,\ref{Sl}, 
and we get, for all $p\in \N$,
\be\label{NIEQnew}\tm U^m\tm_{(t_p,\rho_p)}\leq A(\rho_p,t_p),\  \forall m\in \N, \ \ \mbox{with }\lim_{p\to \infty} A(\rho_p, t_p)=0.\ee   We set $w^{m}:=U^{m}-U^{m-1}$. Hence from \rf{IEQ} we arrive at ($m\geq0$ and $U^{-1}=0$)
$$w^{m+1}(t,x)=-\n_x E*(w^m\otimes U^{m})(t,x)-\n_x E*(U^{m-1}\otimes w^{m})(t,x)\,.$$ Employing the arguments of Lemma\,\ref{LII}, Lemma\,\ref{KRHOH} and Lemma\,\ref{KWHOLE}, and recalling estimate \rf{NIEQ}, we easily arrive at the sequence of estimates 
\be\label{EC}\tm w^1\tm_{(t,\rho)} \leq c_1A^2(\rho,t), \dots, \tm w^m\tm_{(t,\rho)}\leq 2^{m-1}c_1^m A^{m+1}(\rho,t),\dots\,.\ee Since \rf{EIEQ} furnishes $A(\rho,t)<1/2c_1<1$ for all  $t\in (0,T(U_0))$, we get the convergence of the sequence $\{U^m\}$ with respect to the functional $\tm\hskip0.05cm\cdot\hskip0.05cm\tm_{(t,\rho)}$. The uniform  convergence of the sequence of  continuous functions $\{U^m\}$ on $(\eta,T)\times\R^3$  ensures that the limit is a continuous function in $(t,x)\in C((\eta,T)\times\R^3)$. We denote by $U$ the limit. \par
Recalling the definition of  the functional $\tm\hskip0.05cm\cdot\hskip0.05cm\tm_{(t,\rho)}$, by virtue of estimate \rf{NIEQnew}, we deduce 
$$\mbox{for all }\ t_p\to 0,\ \lim_{p\to \infty}\sup_{(0,t_p)}\tau^\frac 12\dm U(\tau)\dm_{\infty}\leq \lim_{p\to \infty}A(\rho_p, t_p)=0,$$
that is estimate \rf{LP}. \par
Further, again from the definition of  the functional $\tm\hskip0.05cm\cdot\hskip0.05cm\tm_{(t,\rho)}$ and using estimate \rf{CRE}, we have, for all $t\in (0, T(U_0)),$
$$\frac{t^\frac 12}{\rho}\sup_{(0,t) } \dm U(\tau)\dm_3\leq 2(h_0+1) \dm U_0\dm_{3,\rho}+ 2(h_1e^{-\frac{\rho^2}{8t}} +\frac{t^\frac 12}{\rho}) \dm U_0\dm_3.$$
Dividing by $t^\frac 12/\rho$ and passing to the limit for $t\to T^-$, we get
$$\sup_{(0,T) } \dm U(\tau)\dm_3\leq \frac{\rho}{T^\frac 12}2(h_0+1) \dm U_0\dm_{3,\rho}+ 2(h_1e^{-\frac{\rho^2}{8T}} \frac{\rho}{T^\frac 12}+1) \dm U_0\dm_3.$$
Hence
$$ \dm U(t)\dm_3\leq \frac{\rho}{T^\frac 12}2(h_0+1) \dm U_0\dm_{3,\rho}+ c \dm U_0\dm_3, \mbox{ for all } t\in (0, T),$$
from which, using that, for any $\rho$, $\dm U_0\dm_{3,\rho}\leq c\dm U_0\dm_3$,  we deduce \eqref{LP1}.
\chiu
 \section{Proof of the main results}In the following, by virtue of estimate \rf{CRE}, we consider 
$A(\rho,t)$ defined in \rf{NIEQ} as a majorant of $\tm U\tm_{(t,\rho)}$, that is \be\label{TMF}\tm U\tm_{(t,\rho)}\leq A(\rho,t)\,,\mbox{ for all }t\in(0,T(U_0))\,,\ee hence, using H\"older's inequality, we get
\be\label{TMF-I}\ba{l}t\dm U(t)\otimes U(t)\dm_\infty\leq A^2(\rho,t)\,,\mbox{ for all }t\in (0,T(U_0))\,,\VS \dm U(t)\otimes U(t)\dm_{\frac32,\rho}+\frac t{\rho^2\hskip-0.1cm}\hskip0.1cm\dm U(t)\otimes U(t)\dm_\frac32\leq A^2(\rho,t)\,,\mbox{ for all }t\in (0,T(U_0))\,.\ea\ee 
{\bf Proof of Theorem\,\ref{EKC}-}   
In the hypothesis of Theorem\,\ref{EKC}, by virtue of Lemma\,\ref{LIEQ} and Lemma\,\ref{CONV}, we establish a solution $U(t,x)$ divergence free to the integral equation \rf{OL} such that for all $t\in[0,T(U_0))$,  $U(t,x)\in J^3(\R^3)$  and satisfies inequality 
\rf{EKC-I}$_3$.  Thanks to Lemma \ref{HHST}, $U$ satisfies the H\"older properties with
\be\label{HPU}t^{\frac12+\frac\theta2}[U(t)]_{\theta}\leq c (A(\rho,t)+A^2(\rho,t))\,,\mbox{ for all }t\in (0,T(U_0))\,.\ee Hence, the following holds:
\be\label{HPUU}t^{1+\frac\theta2}[U(t)\otimes U(t)]_{\theta}\leq c(A^2(\rho,t)+A^3(\rho,t))\,,\mbox{ for all }t\in(0,T(U_0))\,.\ee
As well as, since    $$
 \n U(t,x)\dy =\n H*U_0(t,x)+\lim_{\vep\to 0}\n W_{\vep}(t,x)= \n H*U_0(t,x)+\n W(t,x)\,,$$
applying Lemma\,\ref{LI}    and Lemma\,\ref{KG}, where we mean $w=U\otimes U$,
we arrive at 
\be\label{nU}\displ t\dm \n U(t)\|_\infty\hskip-0.05cm + t^{1+\frac\theta2}[\n U(t)]_{\theta}\hskip-0.05cm\leq \displ \hskip-0.05cm  c A(\rho,t) +c\Big[A^2(\rho,t)+A^3(\rho,t)\Big]\!\leq \hskip-0.05cmc \Big[ A(\rho,t) +A^3(\rho,t)\Big],
\ee where  we employed  \rf{TMF-I}$_2$ and \rf{HPUU}.
Since, from property \texttt{P}, $A(\rho,t)$ tends to zero as $t\to 0$, 
we get estimate \rf{EKC-L}$_1$ for the $\n U$. \par Then, we consider $\pi$  solution to the Poisson equation 
$\Delta \pi=-\n \cdot \n \cdot (U\otimes U)$. We obtain estimates \rf{EKC-L}$_{1,2}$ by applying Lemma\,\ref{LPE-L}. \par 
Since $U$ is solution to the integral equation \rf{OL}, by the couple $(U,\pi)$  one finds the wanted solution to system \rf{NSC}. Concerning the initial condition $U_0$, we firstly observe that the limit property \rf{EKC-L}$_4$ trivially  holds  for $U^0(t,x)$. Then, via the integral equation \rf{OL},  
    and Lemma\,\ref{KWHOLE} for $\n E*(U\otimes U)$, we get
$$\dm U(t)-U^0(t)\dm_{3}\leq c\sup_{(0,t)}\tau^\frac12\dm U(\tau)\dm_\infty\dm U(\tau)\dm_{3}\,\mbox{ for all }t\in[0,T(U_0))$$ Thus, from \rf{LP} and \rf{LP1}, we arrive at the limit property \rf{EKC-L}. \par
Finally,  if we require $\dm U_0\dm_{3}$ sufficiently small and consider $t^\frac12/\rho$ in constant ratio, since $\dm U_0\dm_{3,\rho}\leq \dm U_0\dm_{3}$ for all $\rho>0$, we can satisfy \rf{EIEQ} for arbitrary $\rho$ and then arbitrary $t$. This gives the stated global existence and completes the proof.\chiu 
By virtue of Theorem\,\ref{EKC}, we get\be\label{CSLT} t^\frac12\dm U\otimes U\dm_3\leq t^\frac12\dm U(t)\dm_\infty\dm U(t)\dm_3\leq c\dm U_0\dm_3 A(\rho,t)\,\ee for all $t\in (0,T(U_0))$.\vskip0.1cm
{\bf Proof of Corollary\,\ref{LQ}-}
From Theorem\,\ref{EKC}, the solution $U$ satisfies inequality \rf{SLQ}$_1$ by interpolation. We get \rf{SLQ}$_{2,3}$ by interpolation too. Actually, it suffices to show that $t^\frac12\n U$, $tD^2 U$ and $tU_t$ have $L^3$-norm bounded. For this task, we recall the well known estimates:\be\label{GG}
\dm \n^k H*U_0(t)\dm_3\leq c(k)t^{-\frac k2}\dm U_0\dm_3\,,\ k\in \N.\ee  Then, for $\n W$, defined in \rf{VNW}, we employ \rf{GLT}. In our case $w=U\otimes U$ and $\n\cdot w=U\cdot\n U\,.$ Hence,  we arrive at 
$$\ba{ll}\dy t^\frac12\dm \n W(t)\dm_3  \hskip-0.2cm&\leq c \displ\Big[\sup_{(0,\frac t2)}\tau^\frac 12\dm U(\tau)\dm_{\infty}\dm U(\tau)\dm_3+ \sup_{(\frac t2,t)}\tau\dm \n U(\tau)\dm_{\infty}\dm U(\tau)\dm_3\Big]\VSE \leq c\dm U_0\dm_3(A(\rho,t)+A^3(\rho,t)),\ea 
$$
that, in turn, employing \rf{GG}, implies
\be\label{nUlq} \dm \n U(t)\dm_3\leq ct^{-\frac 12} \dm U_0\dm_3\Big[1+A^3(\rho,t)\Big]. \ee 
By interpolating such estimate and estimate \rf{nU},  we get \rf{SLQ}$_2$. \par Since $\n\cdot U\otimes U=U\cdot\n U$, via the H\"older properties of $U$ and $\n U$, \rf{HPU} and \rf{nU} respectively, we get $U\cdot \n U\in C^{0,\theta}(\R^3)$. Set $\n\cdot w=U\cdot\n U$, via \rf{TMF} and \rf{nU}, since the assumptions of Lemma\,\ref{KG} are satisfied, we apply \rf{GGW} and we deduce  
 $$\ba{ll}t^\frac 32\dm \n\n W_\vep (t)\dm_\infty\hskip-0.2cm& \leq \displ c\sup_{(\frac t2,t)}\tau^{\frac32+\frac\theta2}[U(\tau)\cdot\n U(\tau)]_\theta+c\sup_{(0,\frac t2)}\hskip-0.05cm\hskip-0.05cm\Big[ \dm U(\tau)\dm_{3,\rho}^2+\frac t{\rho^2\hskip-0.1cm}\hskip0.1cm\dm U(\tau)\dm_3^2\Big]\VSE\leq c \Big[A^2(\rho,t)+A^4(\rho,t)\Big].\ea$$
 From this estimate and estimate \rf{LI-I} on $\n\n H*U_0(t,x)$ , we finally obtain  
\be\label{nnU} t^\frac 32\dm \n\n U(t)\|_\infty \leq    c\Big[A(\rho,t)+A^4(\rho,t)\Big],
\ee 
that ensures $$\lim_{t\to 0} t^\frac 32\dm \n\n U(t)\|_\infty  =0.$$
In order to prove that $\n\n W(t)\in L^3(\R^3)$, we apply \rf{GGLT}, taking into account that $\n\n\cdot w=\n (U\cdot\n U)=\n U\cdot\n U+U\cdot\n\n U$. Hence, by means of a  trivial computation,  we find
$$\ba{l}t\displ \dm \n\n W(t)\dm_3 \leq c\Big[\sup_{(0,\frac t2)}\tau^\frac 12\dm U(\tau)\dm_{\infty}\dm U(\tau)\dm_3+ \sup_{(\frac t2,t)}\tau\dm \n U(\tau)\dm_{\infty}\sup_{(\frac t2,t)}\tau^\frac 12\dm \n U(\tau)\dm_3\VS\hskip6cm+\sup_{(\frac t2,t)}\tau^\frac 32\dm \n\n U(\tau)\dm_{\infty}  \dm U(\tau)\dm_3\Big]. \ea
$$
Hence, by applying estimates  \rf{TMF}, \rf{nU}, \rf{nUlq} and \rf{nnU}, we get
$$ t\dm \n\n W(t)\dm_3\leq c \dm U_0\dm_3\Big[A(\rho,t)+A^6(\rho,t)\Big], $$
that, in turn, employing \rf{GG}, implies
$$ t\dm \n\n U(t)\dm_3\leq c \dm U_0\dm_3\Big[1+A^6(\rho,t)\Big]. $$
By interpolating such estimate and estimate \rf{nnU},  we get \rf{SLQ}$_3$. \par

\chiu

{\bf{Proof of Theorem\,\ref{PU} }} Let us consider two solutions $U$ and $\ov U$ satisfying te assumptions, and set $u:=\ov U-U$ . Then, for all $t>s\geq 0$, $u$ satisfies the following integral equation
\be\label{VW}\ba{ll}\displaystyle\intl{s}^t(u(\tau),\,\varphi_\tau+\Delta \varphi)\,d\tau +\intl{s}^t\Big[(u(\tau)\cdot\nabla \varphi, \ov U)+(U\cdot \nabla\varphi, u)\Big]\,d\tau\VS\hfill =(u(t), \varphi(t))-(u(s), \varphi(s)), \ \forall \varphi\in C^1([0,T]; \mathscr C_0(\R^3)).\ea \ee
We denote by $\psi$ the solution to the Cauchy problem:
\be\label{STE}\ba{l}\psi_t-\Delta \psi=-\n\pi_\psi\,,\quad \n\cdot \psi=0\,\mbox{ in }(0,T)\times\R^3\,,\VS\psi=\psi_0\in\mathscr C_0(\R^3)\mbox{ on }\{0\}\times\R^3\,.\ea\ee 
  It is well known that $\psi$ is a smooth solution with $\psi\in C([0,T);J^p(\R^3))$, for all $p\in(1,\infty)$, and satisfies the following estimates:\be\label{H-I}\ba{ll}\hskip0.2cmq\geq p\,,\,\dm \psi(t)\dm_q&\leq ct^{-\frac32\left(\frac1p-\frac1q\right)}\dm \psi_0\dm_p\,,\;\mbox{ for all }t>0\,,\VS \null\hskip1cm\dm\n \psi(t)\dm_q&\leq c_1t^{-\frac12-\frac32\left(\frac1p-\frac1q\right)}\dm\psi_0\dm_p\,,\mbox{ for all }t>0\,.\ea\ee  For   $t>0$, we set $\widehat\psi(\tau,x):=\psi(t-\tau,x)$ provided that $(\tau,x)\in(0,t)\times\R^3$. It is well known that $\widehat\psi$ is a solution backward in time with $\widehat\psi(t,x)=\psi_0(x)$. 
  \par 
  Let us write the integral equation \rf{VW} with $\widehat\psi$ in place of $\varphi$. We get 
\be\label{ADR}(u(t),\psi_0)=(u(s),\widehat\psi
(s))+\intll
st\Big[ (U\cdot\nabla\widehat\psi,u)+( u\cdot \n \widehat\psi,\ov U)\Big]d\tau .\ee
 Hence
 \be\label{UES}\hskip-0.15cm\ba{ll}|(u(t),\psi_0)|\displ \displ \leq |(u(s),\widehat\psi
(s))|\VS\hskip1.75cm+\!  c\sup_{(s,t)}\Big[\tau^\frac12(\dm U\big(\tau)\dm_\infty\!+\dm \ov U(\tau)\dm_\infty\big)\Big]\sup_{(s,t)}\dm u(\tau)\dm_{3}
 \intll {s}{t}\!\! \tau^{-\frac 12}\dm \n \psi(t\!-\!\tau)\dm_\frac 32d\tau\VS\displaystyle \hskip1.7cm\leq 
 |(u(s),\widehat\psi
(s))|\VS\hskip1.75cm+
 c \dm \psi_0\dm_\frac 32\sup_{(s, t)}\tau^\frac 12(\dm U(\tau)\dm_\infty\!+ \dm \ov U(\tau)\dm_\infty)\sup_{(s,t)}\dm u(\tau)\dm_{3}\! 
  \intll {s}{t}\!\!\tau^{-\frac 12}(t\!-\!\tau)^{-\frac12}d\tau,  
\ea\ee   
  for all $t\in [0,T(U))\cap[0,T(\ov U)) $. Since $\psi_0$ is arbitrary, and then letting  $s\to 0$, we obtain
  $$\dm u(t)\dm_3\leq  c \dm \sup_{(0, t)}\Big[\tau^\frac 12(\dm U(\tau)\dm_\infty+ \dm \ov U(\tau)\dm_\infty)\Big]\sup_{(0,t)}\dm u(\tau)\dm_{3} .
$$   
  From the validity of the limit property \rf{LPB}$_1$ on both solutions, one easily deduces the uniqueness on some interval $(0,\delta]$. It remains to discuss the uniqueness when $t\geq \delta$. Writing estimate \rf{UES} with $s=\delta$, since $\dm u(\delta)\dm_3=0$, 
 we deduce the estimate $$\ba{ll}\dm u(t)\dm_3&\displ\hskip-0.05cm\leq \hskip-0.05cm\!\!\intll {\delta}{t}\!\! \Big[\dm U(\tau)\dm_\infty+\dm \ov U(\tau)\dm_\infty\Big]\hskip-0.05cm\dm u(\tau)\dm_3(t-\tau)^{-\frac12}d\tau\\&\hskip-0.05cm\leq  c\, \displ \delta^{-\frac 12}\sup_{(\delta, t)}\tau^\frac 12(\dm U(\tau)\dm_\infty+ \dm \ov U(\tau)\dm_\infty) \intll {\delta}{t}\dm u(\tau)\dm_3(t-\tau)^{-\frac12}\,d\tau\,.  
\ea$$   We are in the hypothesis of the logarithmic Gronwall inequality (see Lemma 4 in \cite{giga}). Therefore we obtain $\dm u(t)\dm_3=0$, for any $t\in [\delta, T)$, that completes the proof. \chiu

\end{document}